\documentclass{colt2017} 

\usepackage{pdfsync,nicefrac}
\usepackage{float}
\usepackage{autonum}
\usepackage{amsfonts,mathrsfs}
\usepackage{dsfont}
\usepackage{cleveref}

\def\ds{\displaystyle}

\newcommand{\RR}{\mathbb{R}}

\newcommand{\NN}{\mathbb{N}}

\newcommand{\bzeta}{\boldsymbol{\zeta}}
\newcommand{\btheta}{\boldsymbol{\theta}}
\newcommand{\bvartheta}{\boldsymbol{\vartheta}}

\newcommand{\bDelta}{\boldsymbol{\Delta}}

\newcommand{\bxi}{\boldsymbol{\xi}}

\newcommand{\bfE}{\mathbf E}

\newcommand{\bu}{\boldsymbol u}

\newcommand{\bU}{\boldsymbol U}
\newcommand{\bV}{\boldsymbol V}
\newcommand{\bW}{\boldsymbol W}

\newcommand{\bL}{\boldsymbol L}

\newcommand{\bx}{\boldsymbol x}
\newcommand{\bY}{\boldsymbol Y}

\newtheorem{lem}{Lemma}

\title[Further analogy between sampling and optimization]
{Further and stronger analogy between sampling and optimization: Langevin Monte Carlo and gradient descent}
\usepackage{times}

 \coltauthor{\Name{Arnak S. Dalalyan} \Email{arnak.dalalyan@ensae.fr}\\
 \addr ENSAE/CREST/Universit\'e Paris Saclay}

\begin{document}

\maketitle

\begin{abstract}
In this paper\footnote{This paper has been published in proceedings of COLT 2017. However, this version 
is more recent. We have corrected some typos ($2/(m+M)$ instead of $1/(m+M)$ on pages 3-4) and slightly
improved the upper bound of \Cref{thTwo}.}, we revisit the recently established theoretical guarantees 
for the convergence of the Langevin Monte Carlo algorithm of sampling from 
a smooth and (strongly) log-concave density. We improve the existing results
when the convergence is measured in the Wasserstein distance and provide 
further insights on the very tight relations between, on the one hand, the 
Langevin Monte Carlo for sampling and, on the other hand, the gradient 
descent for optimization. Finally, we also establish guarantees for the
convergence of a version of the Langevin Monte Carlo algorithm that is based
on noisy evaluations of the gradient. 
\end{abstract}


\section{Introduction}

Let $p$ be a positive integer and $f:\RR^p\to\RR$ be a measurable function such that the 
integral $\int_{\RR^p} \exp\{-f(\btheta)\}\,d\btheta$ is finite. In various applications, 
one is faced with the problems of finding the minimum point of $f$ or computing the average 
with respect to the probability density 
\begin{equation}
\pi(\btheta) = \frac{e^{-f(\btheta)}}{\int_{\RR^p} e^{-f(\bu)}\,d\bu}.
\end{equation}
In other words, one often looks for approximating the values $\btheta^*$ and $\bar\btheta$
defined as
\begin{equation} 
\bar\btheta = \int_{\RR^p} \btheta\,\pi(\btheta)\,d\btheta,\qquad
\btheta^*\in\text{arg}\min_{\btheta\in\RR^p} f(\btheta).
\end{equation}
In most situations, the approximations of these values are computed using iterative algorithms 
which share many common features. There is a vast variety of such algorithms for solving both tasks, 
see for example \citep{BoydBook} for optimization and \citep{Atchade2011} for approximate sampling. 
The similarities between the task of optimization and that of averaging have been recently exploited
in the papers \citep{Dalalyan14,Durmus2,Durmus1} in order to establish fast and accurate theoretical
guarantees for sampling from and averaging with respect to the density $\pi$ using the Langevin
Monte Carlo algorithm. The goal of the present work is to push further this study both by improving 
the existing bounds and by extending them in some directions.

We will focus on strongly convex functions $f$ having a Lipschitz continuous gradient. 
That is, we assume that there exist two positive constants $m$ and $M$ such that 
\begin{equation} \label{1}
\begin{cases}
f(\btheta)-f(\btheta')-\nabla f(\btheta')^\top (\btheta-\btheta')
\ge (\nicefrac{m}2)\|\btheta-\btheta'\|_2^2, \text{\vphantom{$I_{\textstyle\int_{I_I}}$}}\\
\|\nabla f(\btheta)-\nabla f(\btheta')\|_2 \le M \|\btheta-\btheta'\|_2,
\end{cases}
\qquad \forall \btheta,\btheta'\in\RR^p,
\end{equation}
where $\nabla f$ stands for the gradient of $f$ and $\|\cdot\|_2$ is the Euclidean norm. 
We say that the density $\pi(\btheta)\propto e^{-f(\btheta)}$ is log-concave (resp.\ strongly 
log-concave) if the function $f$ satisfies the first inequality of (\ref{1}) with $m=0$ 
(resp.\ $m>0$).

The Langevin Monte Carlo (LMC) algorithm studied throughout this work is the analogue of the 
gradient descent algorithm for optimization. Starting from an initial point $\bvartheta^{(0)}\in\RR^p$ 
that may be deterministic or random, the iterations of the algorithm are defined by
the update rule
\begin{align}\label{2}
\bvartheta^{(k+1,h)} = \bvartheta^{(k,h)} - h \nabla f(\bvartheta^{(k,h)})+ \sqrt{2h}\;\bxi^{(k+1)};
\qquad k=0,1,2,\ldots
\end{align}
where $h>0$ is a tuning parameter, referred to as the step-size, and $\bxi^{(1)},\ldots,\bxi^{(k)},\ldots$
is a sequence of mutually independent, and independent of $\bvartheta^{(0)}$, centered Gaussian vectors with
covariance matrices equal to identity. 
Under the assumptions imposed on $f$, when $h$ is small and $k$ is large (so that the product 
$kh$ is large), the distribution of $\bvartheta^{(k,h)}$ is close in various metrics to the 
distribution with density $\pi(\btheta)$, hereafter referred to as the target distribution. 
An important question is to quantify this closeness; this might be particularly useful for
deriving a stopping rule for the LMC algorithm.

The measure of approximation used in this paper is the Wasserstein-Monge-Kantorovich 
distance $W_2$. For two measures $\mu$ and $\nu$ defined on $(\RR^p,\mathscr B(\RR^p))$, 
$W_2$ is defined by
\begin{equation}
W_2(\mu,\nu) = \Big(\inf_{\gamma\in \Gamma(\mu,\nu)} \int_{\RR^p\times \RR^p} 
		\|\btheta-\btheta'\|_2^2\,d\gamma(\btheta,\btheta')\Big)^{1/2},
\end{equation}
where the $\inf$ is with respect to all joint distributions $\gamma$ having $\mu$ and
$\nu$ as marginal distributions. This distance is perhaps more suitable for quantifying
the quality of approximate sampling schemes than other metrics such as the total variation. 
Indeed, on the one hand, bounds on the Wasserstein distance---unlike the bounds on
the total-variation distance---directly provide the level of approximating 
the first order moment. For instance, if $\mu$ and $\nu$ are two Dirac measures at 
the points $\btheta$ and $\btheta'$, respectively, then the total-variation distance
$D_{\rm TV}(\delta_{\btheta},\delta_{\btheta'})$ equals one whenever $\btheta\not=\btheta'$,
whereas $W_2(\delta_{\btheta},\delta_{\btheta'}) = \|\btheta-\btheta'\|_2$ is a smoothly
increasing function of the Euclidean distance between $\btheta$ and $\btheta'$. This
seems to better correspond to the intuition on the closeness of two distributions.

\section{Improved guarantees for the Wasserstein distance}

The rationale behind the LMC algorithm \eqref{2} is simple: the Markov chain 
$\{\bvartheta^{(k,h)}\}_{k\in\NN}$ is the Euler discretization of a continuous-time
diffusion process $\{\bL_t :t\in\RR_+\}$, known as Langevin diffusion, that has $\pi$ 
as invariant density  \citep[Thm. 3.5]{bhattacharya1978}. The Langevin diffusion is 
defined by the stochastic differential equation 
\begin{equation}\label{3}
d\bL_t = -\nabla f(\bL_t)\,dt + \sqrt{2} \; d\bW_t,\qquad t\ge 0,
\end{equation}
where $\{\bW_t:t\ge 0\}$ is a $p$-dimensional Brownian motion. When $f$ satisfies 
condition (\ref{1}), equation (\ref{3}) has a unique
strong solution which is a Markov process. Let $\nu_k$ be the distribution of the
$k$-th iterate of the LMC algorithm, that is $\vartheta^{(k,h)}\sim \nu_k$.

\begin{theorem}\label{thOne}
Assume that $h\in(0,\nicefrac2M)$. The following claims hold:
\begin{enumerate}\itemsep=10pt
\item[{\rm (a)}] If $h\le \nicefrac2{(m+M)}$ then $W_2(\nu_K, \pi) \le 
(1-mh)^K W_2(\nu_0,\pi) + 1.82(M/m)(hp)^{1/2}$.
\item[{\rm (b)}] If $h\ge \nicefrac2{(m+M)}$ then $W_2(\nu_K, \pi) \le 
\ds (Mh-1)^K W_2(\nu_0,\pi) + 1.82\frac{Mh}{2-Mh}(hp)^{1/2}$.
\end{enumerate}
\end{theorem}

The proof of this theorem is postponed to \Cref{secProof}. We content ourselves here 
by discussing the relation of this result to previous work. Note that if the initial 
value $\bvartheta^{(0)}=\btheta^{(0)}$ is deterministic then, according to 
\citep[Theorem 1]{Durmus2}, we have
\begin{align}
W_2(\nu_0,\pi)^2 
	& = \int_{\RR^p} \|\btheta^{(0)}-\btheta\|_2^2\pi(d\btheta)\\
	& = \|\btheta^{(0)}-\bar\btheta\|_2^2 + \int_{\RR^p} \|\bar\btheta-\btheta\|_2^2\pi(d\btheta)\\
	& \le \|\btheta^{(0)}-\bar\btheta\|_2^2 + p/m.\label{4}
\end{align}

First of all, let us remark that if we choose $h$ and $K$ so that 
\begin{equation}\label{5}
h\le \nicefrac{2}{(m+M)},\qquad e^{-mhK}W_2(\nu_0,\pi)\le \varepsilon/2,\quad 
1.82(M/m)(hp)^{1/2}\le \varepsilon/2,
\end{equation}
then we have $W_2(\nu_K, \pi) \le \varepsilon$. In other words, conditions
\eqref{5} are sufficient for the density of the output of the LMC algorithm with $K$ 
iterations to be within the precision $\varepsilon$ of the target density when the precision
is measured using the Wasserstein distance. This readily yields 
\begin{equation}\label{6}
h\le \frac{m^2\varepsilon^2}{14M^2p}\wedge \frac2{m+M}\quad\text{and}\quad 
hK\ge \frac1m\log\Big(\frac{2(\|\btheta^{(0)}-\bar\btheta\|_2^2+ p/m)^{1/2}}\varepsilon\Big)
\end{equation}
Assuming $m,M$ and $\|\btheta^{(0)}-\bar\btheta\|_2^2/p$ to be constants, we can deduce
from the last display that it suffices $K = C p\varepsilon^{-2}\log(p/\varepsilon)$
number of iterations in order to reach the precision level $\varepsilon$. This fact has
been first established in \citep{Dalalyan14} for the LMC algorithm with a warm start
and the total-variation distance. It was later improved by \cite{Durmus2}, who showed
that the same result holds for any starting point and established similar bounds for
the Wasserstein distance.

In order to make the comparison easier, let us recall below the corresponding result 
from\footnote{We slightly adapt the original result taking into account the fact
that we are dealing with the LMC algorithm with a constant step.}
\citep{Durmus2}. It asserts that under condition \eqref{1}, if $h\le \nicefrac2{(m+M)}$ 
then
\begin{equation}\label{l}
W_2^2(\nu_{K}, \pi) \le 
2\Big(1-\frac{mMh}{m+M}\Big)^K W^2_2(\nu,\pi) + 
\frac{Mhp}{m}(m+M)\Big(h + \frac{m+M}{2mM}\Big)\Big(2+\frac{M^2h}{m}+\frac{M^2h^2}{6}\Big).
\end{equation}
When we compare this inequality with the claims of \Cref{thOne}, we see that 
\begin{enumerate}
\item[i)] \Cref{thOne} holds under weaker conditions: $h\le \nicefrac2{M}$ instead of 
$h\le \nicefrac2{(m+M)}$. 
\item[ii)] The analytical expressions of the upper bounds on the Wasserstein distance 
in \Cref{thOne} are not as involved as those of \eqref{l}.
\item[iii)] If we take a closer look, we can check that when  $h\le \nicefrac2{(m+M)}$, 
the upper bound in part (a) of \Cref{thOne} is sharper than that of \eqref{l}.
\end{enumerate} 

\begin{figure}
\centerline{\includegraphics[width = 0.75\textwidth]{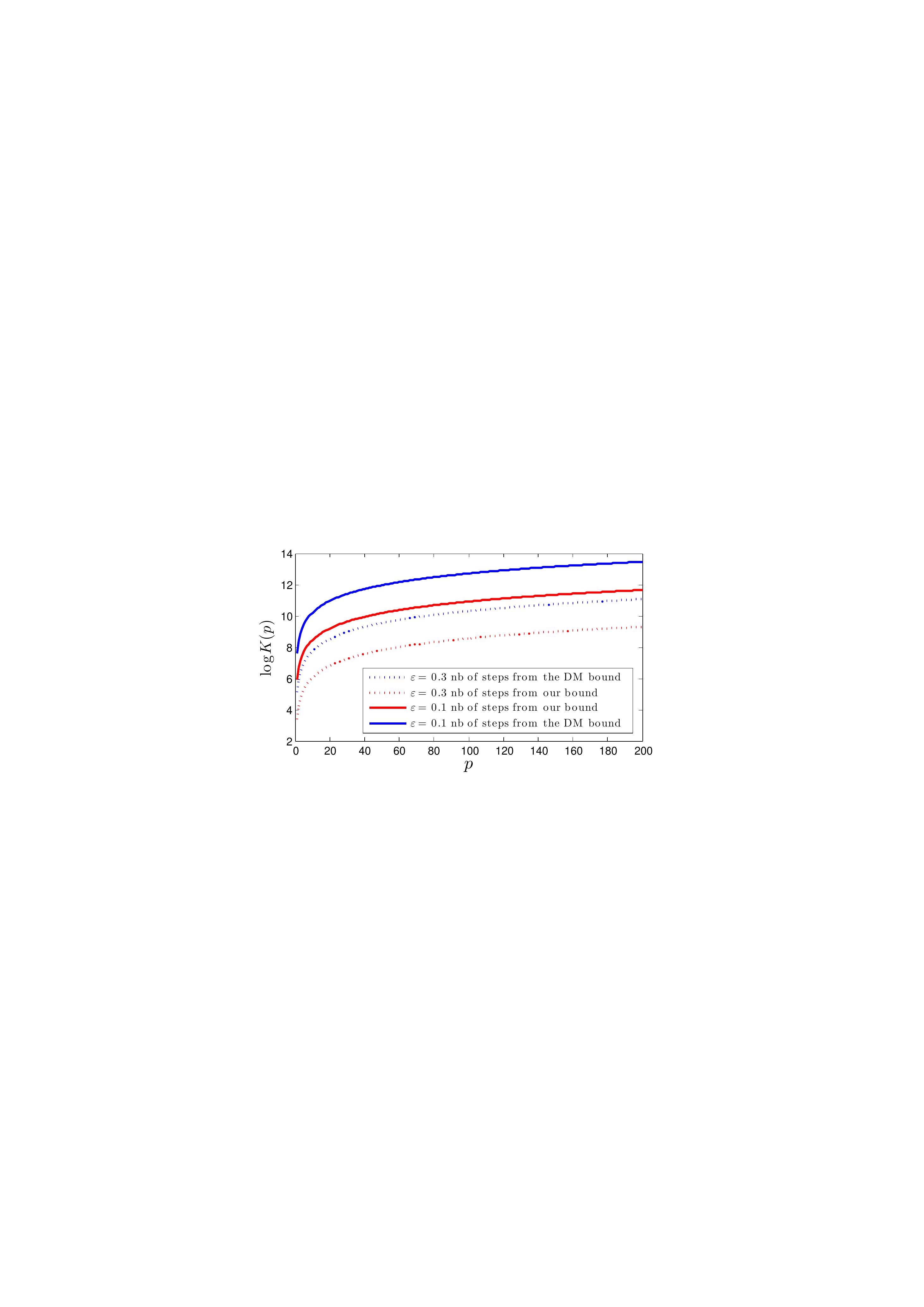}}
\caption{The curves of the functions $p\mapsto \log K(p)$, where $K(p)$ is the number of steps--- 
derived either from our bound or from the bound \eqref{l} of \citep{Durmus2}---sufficing for 
reaching the precision level $\varepsilon$ (for $\varepsilon = 0.1$ and $\varepsilon = 0.3$).}
\label{figOne}
\end{figure}

In order to better illustrate the claim in iii) above, we consider a numerical example
in which $m=4$, $M = 5$ and $\|\btheta^{(0)}-\bar\btheta\|_2^2 = p$. 
Let $F_{\rm our}(h,K,p)$ and $F_{\rm DM}(h,K,p)$ be the upper bounds on $W_2(\nu_K,\pi)$ 
provided by \Cref{thOne} and \eqref{l}. For different values of $p$, we compute
\begin{align}
K_{\rm our}(p) &= \min \big\{K  : \text{ there exists $h\le \nicefrac2{(m+M)}$ such that } 
								F_{\rm our}(h,K,p)\le \varepsilon\big\},\\
K_{\rm DM}(p) &= \min \big\{K  : \text{ there exists $h\le \nicefrac2{(m+M)}$ such that } 
								F_{\rm DM}(h,K,p)\le \varepsilon\big\}.								
\end{align} 
The curves of the functions $p\mapsto \log K_{\rm our}(p)$ and $p\mapsto \log K_{\rm DM}(p)$, 
for $\varepsilon = 0.1$ and $\varepsilon = 0.3$ are plotted in~\Cref{figOne}. We can 
deduce from these plots that the number of iterations yielded by our bound is more than 5 
times smaller than the number of iterations recommended by bound  \eqref{l} of 
\cite{Durmus2}.

\begin{remark}
Although the upper bound on $W_2(\nu_0,\pi)$ provided by~\eqref{4} is relevant for
understanding the order of magnitude of  $W_2(\nu_0,\pi)$, it has limited applicability 
since the distance $\|\btheta_0-\bar\btheta\|$ might be hard to evaluate. An attractive 
alternative to that bound is the following\footnote{The second line follows from
strong convexity whereas the third line is a consequence of the two identities 
$\int_{\RR^p} \nabla f(\btheta)\pi(d\btheta) = 0$ and $\int_{\RR^p} \btheta^\top
\nabla f(\btheta)\pi(d\btheta) = p$. These identities follow from the fundamental 
theorem of calculus and the integration by parts formula, respectively.}:
\begin{align}
W_2(\nu_0,\pi)^2 
	& = \int_{\RR^p} \|\btheta^{(0)}-\btheta\|_2^2\pi(d\btheta)\\
	& \le \frac2{m}\int_{\RR^p}\Big(f(\btheta_0)-f(\btheta)-\nabla f(\btheta)^\top(\btheta_0-\btheta)
	\Big)\pi(d\btheta)\\
	& = \frac2{m}\Big(f(\btheta_0)-\int_{\RR^p} f(\btheta)\,\pi(d\btheta)+p\Big).\label{init1}
\end{align}
If $f$ is lower bounded by some known constant, for instance if $f\ge 0$, the last inequality
provides the computable upper bound $W_2(\nu_0,\pi)^2 \le  \frac2{m}\big(f(\btheta_0)+p\big)$.
\end{remark}

\section{Relation with optimization}
\label{secOpt}

We have already mentioned that the LMC algorithm is very close to the gradient descent 
algorithm for computing the minimum $\btheta^*$ of the function $f$.  However, when we 
compare the guarantees of~\Cref{thOne} with those available for the optimization problem, 
we remark the following striking difference. The approximate computation of $\btheta^*$ requires
a number of steps of the order of $\log(1/\varepsilon)$ to reach the precision $\varepsilon$,
whereas, for reaching the same precision in sampling from $\pi$, the LMC algorithm needs a 
number of iterations proportional to $(p/\varepsilon^2)\log (p/\varepsilon)$. 
The goal of this section is to explain that this, at first sight very disappointing 
behavior of the LMC algorithm is, in fact, continuously connected to the exponential 
convergence of the gradient descent.

The main ingredient for the explanation is that the function $f(\btheta)$ and the function 
$f_\tau(\btheta) = f(\btheta)/\tau$ have the same point of minimum $\btheta^*$, whatever 
the real number $\tau>0$. In addition, if we define the density function 
$\pi_\tau(\btheta)\propto \exp\big(-f_\tau(\btheta)\big)$, then the average value
$$
\bar\btheta_\tau = \int_{\RR^p } \btheta\, \pi_\tau(\btheta)\,d\btheta
$$
tends to the minimum point $\btheta^*$ when $\tau$ goes to zero. Furthermore, 
the distribution $\pi_\tau(d\btheta)$ tends to the Dirac measure at $\btheta^*$.
Clearly, $f_\tau$ satisfies \eqref{1} with the constants $m_\tau = m/\tau$ 
and $M_\tau = M/\tau$. Therefore, on the one hand, we can apply to $\pi_\tau$ 
claim (a) of \Cref{thOne}, which tells  us that if we choose $h = 1/M_\tau = \tau/M$, 
then
\begin{equation}
\label{7}
W_2(\nu_K,\pi_\tau) \le \Big(1-\frac{m}{M}\Big)^K W_2(\delta_{\btheta^{(0)}},\pi_\tau)
+ 2\Big(\frac{M}{m}\Big)\Big(\frac{p\tau}{M}\Big)^{1/2}.
\end{equation}
On the other hand, the LMC algorithm with the step-size $h=\tau/M$ applied to 
$f_\tau$ reads as
\begin{equation} 
\label{8}
\bvartheta^{(k+1,h)} = \bvartheta^{(k,h)} - \frac1M \nabla f(\bvartheta^{(k,h)})+ 
\sqrt{\frac{2\tau}M}\;\bxi^{(k+1)};\qquad k=0,1,2,\ldots
\end{equation} 
When the parameter $\tau$ goes to zero, the LMC sequence \eqref{8} tends
to the gradient descent sequence $\btheta^{(k)}$. Therefore, the limiting case
of \eqref{7} corresponding to $\tau\to 0$ writes as 
\begin{equation}
\label{optimGuar}
\|\btheta^{(K)}-\btheta^*\|_2 \le \Big(1-\frac{m}{M}\Big)^K \|\btheta^{(0)}-\btheta^*\|_2,
\end{equation}
which is a well-known result in Optimization. This clearly shows that \Cref{thOne} is a natural
extension of the results of convergence from optimization to sampling. 

\section{Guarantees for the noisy gradient version}

In some situations, the precise evaluation of the gradient $\nabla f(\btheta)$ 
is computationally expensive or practically impossible, but it is possible to
obtain noisy evaluations of $\nabla f$ at any point. This is the setting considered
in the present section. More precisely, we assume that at any point $\bvartheta^{(k,h)}\in\RR^p$
of the LMC algorithm, we can observe the value
\begin{equation}
\bY^{(k,h)}  = \nabla f(\bvartheta^{(k,h)}) + \sigma\,\bzeta^{(k)},
\end{equation} 
where $\{\bzeta^{(k)}:\,k=0,1,\ldots\}$ is a sequence of independent zero mean random vectors 
such that $\bfE[\|\bzeta^{(k)}\|_2^2]\le p$ and $\sigma>0$ is a deterministic noise level. Furthermore, 
the noise vector $\bzeta^{(k)}$ is independent of the past states 
$\bvartheta^{(1,h)},\ldots,\bvartheta^{(k,h)}$. The noisy LMC (nLMC)  algorithm is then defined as
\begin{align}\label{9}
\bvartheta^{(k+1,h)} = \bvartheta^{(k,h)} - h \bY^{(k,h)}+ \sqrt{2h}\;\bxi^{(k+1)};\qquad k=0,1,2,\ldots
\end{align}
where $h>0$ and $\bxi^{(k+1)}$ are as in \eqref{2}. The next theorem extends the guarantees 
of \Cref{thOne} to the noisy-gradient setting and to the nLMC algorithm.

\begin{theorem}\label{thTwo}
Let $\bvartheta^{(K,h)}$ be the $K$-th iterate of the nLMC algorithm \eqref{9} and 
$\nu_K$ be its distribution. If the function $f$ satisfies condition \eqref{1} 
and $h\le 2/M$ then the following claims hold:
\begin{enumerate}\itemsep=2pt
\item[{\rm (a)}] If $h\le \nicefrac2{(m+M)}$ then 
\begin{align}\label{A}
W_2(\nu_K, \pi) \le 
\Big(1-\frac{mh}{2}\Big)^{K} W_2(\nu_0,\pi) + 
		\Big(\frac{2hp}{m}\Big)^{1/2}\Big\{\sigma^2  + \frac{3.3M^2}{m}\Big\}^{1/2}.
\end{align}
\item[{\rm (b)}] If $h\ge \nicefrac2{(m+M)}$ then 
$$
W_2(\nu_K, \pi) 
\le \Big(\frac{Mh}{2}\Big)^{K} W_2(\nu_0,\pi) + 
		\Big(\frac{2h^2p}{2-Mh}\Big)^{1/2}\Big\{\sigma^2  + \frac{6.6M}{2-Mh}\Big\}^{1/2}.
$$
\end{enumerate}
\end{theorem}

To understand the potential scope of applicability of this result, let us consider a typical 
statistical problem in which $f(\btheta)$ is the negative log-likelihood of $n$ independent
random variables $X_1,\ldots,X_n$. Then, if $\ell(\btheta,x)$ is the log-likelihood of one 
variable, we have
$$
f(\btheta) = \sum_{i=1}^n \ell(\btheta,X_i).
$$                     
In such a situation, if the Fisher information is not degenerated, both $m$ and $M$ are
proportional to the sample size $n$. When the gradient of $\ell(\btheta,X_i)$ 
with respect to parameter $\btheta$ is hard to compute, one can replace the evaluation
of $\nabla f(\bvartheta^{(k,h)})$ at each step $k$ by that  of $Y_k= n \nabla_{\btheta} \ell 
(\bvartheta^{(k,h)},X_k)$. Under suitable assumptions, this random vector satisfies
the conditions of \Cref{thTwo} with a $\sigma^2$ proportional to $n$. Therefore, if 
we analyze the expression between curly brackets in \eqref{A}, we see that the
additional term, $\sigma^2$, due to the subsampling is 
of the same order of magnitude as the term $3.3M^2/m$. Thus, using the subsampled 
gradient in the LMC algorithm does not cause a significant deterioration of the
precision while reducing considerably the computational burden.

\section{Discussion and outlook}

We have established simple guarantees for the convergence of the Langevin Monte Carlo 
algorithm under the Wasserstein metric. These guarantees are valid under strong convexity
and Lipschitz-gradient assumptions on the log-density function, for a step-size 
smaller than $2/M$, where $M$ is the constant in the Lipschitz condition. These guarantees
are sharper than previously established analogous results and in perfect agreement with
the analogous results in Optimization. Furthermore, we have shown that similar results
can be obtained in the case where only noisy evaluations of the gradient are possible.

There are a number of interesting directions in which this work can be extended. 
One relevant and closely related problem is the approximate computation of the volume
of a convex body, or, the problem of sampling from the uniform distribution on a 
convex body. This problem has been analyzed by other Monte Carlo methods such as 
``Hit and Run'' in a series of papers by \cite{Lovasz2,Lovasz1}, see also the more 
recent paper \citep{Bubeck15}. Numerical experiments reported in \citep{Bubeck15} suggest 
that the LMC algorithm might perform better in practice than ``Hit and Run''. It would 
be interesting to have a theoretical result corroborating this observation.

Other interesting avenues for future research include the possible adaptation of the 
Nesterov acceleration to the problem of sampling, extensions to second-order methods as
well as the alleviation of the strong-convexity assumptions. We also plan to investigate
in more depth the applications is high-dimensional statistics (see, for instance, 
\cite{DalalyanTsybakov12a}). Some results in these directions are already obtained 
in \citep{Dalalyan14,Durmus2,Durmus1}. It is a stimulating  question whether we can 
combine ideas of the present work and the aforementioned earlier results to get improved 
guarantees.

\section{Proofs}
\label{secProof}

The first part of the proofs of \Cref{thOne} and \Cref{thTwo} is the same. 
We start this section by this common part and then we proceed with the
proofs of the two theorems separately. 

Let $\bW$ be a $p$-dimensional Brownian Motion such that $\bW_{(k+1)h} - \bW_{kh} = \sqrt{h}\,
\bxi^{(k+1)}$. We define the stochastic process $\bL$ 
so that $\bL_0\sim \pi$ and 
\begin{align}\label{B}
\bL_t &= \bL_0 - \int_0^t  \nabla f(\bL_s)\,ds + \sqrt{2}\,\bW_t,\qquad\forall\, t>0.
\end{align}
It is clear that this equation implies that
\begin{align}
\ds\bL_{(k+1)h} 
	&= \bL_{kh} - \int_{kh}^{(k+1)h} \nabla f(\bL_s)\,ds + \sqrt{2}\,(\bW_{(k+1)h}-\bW_{kh})\\
	&= \bL_{kh} - \int_{kh}^{(k+1)h} \nabla f(\bL_s)\,ds + \sqrt{2h}\,\bxi^{(k+1)}.
\end{align}
Furthermore, $\{\bL_t:t\ge 0\}$ is a diffusion process having $\pi$ as the stationary
distribution. Since the initial value $\bL_0$ is drawn from $\pi$, we have $\bL_t\sim \pi$
for every $t\ge 0$.

Let us denote $\bDelta_k = \bL_{kh}-\bvartheta^{(k,h)}$ and $I_k = (kh,(k+1)h]$. We have
\begin{align}
\bDelta_{k+1} 
	& = \bDelta_k  + h \bY^{(k,h)} - \int_{I_k}\nabla f(\bL_t)\,dt \\
	& = \bDelta_k  - h\big(\underbrace{\nabla f(\bvartheta^{(k,h)}+\bDelta_k)-\nabla f(\bvartheta^{(k,h)})
			}_{:=\bU_k}\big)+ \sigma h\bzeta^{(k)}
			-\underbrace{\int_{I_k}\big(\nabla f(\bL_t) - \nabla f(\bL_{kh})\big)\,dt}_{:=\bV_k}.
\end{align}
In view of the triangle inequality, we get
\begin{equation}
\|\bDelta_{k+1} \|_2  \le \|\bDelta_k  -h \bU_k + \sigma h \bzeta^{(k)}\|_2 + \|\bV_k\|_2.\label{C}
\end{equation}
For the first norm in the right hand side, we can use the following inequalities:
\begin{align}
\bfE[\|\bDelta_k  -h \bU_k + \sigma h \bzeta^{(k)}\|_2^2]
		& = \bfE[\|\bDelta_k  -h \bU_k\|_2^2] + \bfE[\|\sigma h \bzeta^{(k)}\|_2^2] \\
		& = \bfE[\|\bDelta_k  -h \bU_k\|_2^2] + \sigma^2 h^2 p.\label{D}
\end{align}

We need now three technical lemmas the proofs of which are postponed to \Cref{ssecLem}. 

\begin{lem}\label{lemA}
Let us introduce the constant $\gamma$ that equals $|1-mh|$ if 
$h\le \nicefrac2{(m+M)}$ and $|1-Mh|$ if $h\ge \nicefrac2{(m+M)}$. (Since $h\in(0,\nicefrac2M)$,
this value $\gamma$ satisfies $0< \gamma <1$). It holds that
\begin{align}
\|\bDelta_k  -h \bU_k\|_2
		& \le  \gamma\|\bDelta_k\|_2.\label{E}
\end{align}
\end{lem}

\begin{lem}\label{lemB}
If the function $f$ is continuously differentiable and the gradient of $f$ is Lipschitz 
with constant $M$, then 
\begin{equation}
\int_{\RR^p} \|\nabla f(\bx)\|_2^2\,\pi(\bx)\,d\bx \le Mp. 
\end{equation}
\end{lem}

\begin{lem}\label{lemC}
If the function $f$ has a Lipschitz-continuous gradient with the Lipschitz constant $M$, 
$\bL$ is the Langevin diffusion \eqref{B} and $\bV(a) = 
\int_a^{a+h}\big(\nabla f(\bL_t)-\nabla f(\bL_a)\big)\,dt$ for some $a\ge 0$, then 
\begin{align}		
\big(\bfE[\|\bV(a)\|^2_2]\big)^{1/2}	&\le \bigg(\frac13 h^4 M^{3}p\bigg)^{1/2} + (h^3p)^{1/2}M . 
\end{align}
\end{lem}
This completes the common part of the proof. We present below the proofs of the theorems.

\subsection{Proof of \Cref{thOne}}

Using \eqref{C} with $\sigma=0$ and \Cref{lemA}, we get
\begin{equation}
\|\bDelta_{k+1} \|_2  \le \gamma \|\bDelta_k\|_2 + \|\bV_k\|_2,\qquad \forall k\in\NN.
\end{equation}
In view of the Minkowski inequality and \Cref{lemC}, this yields 
\begin{align}
(\bfE[\|\bDelta_{k+1} \|_2^2])^{1/2}  
		&\le \gamma (\bfE[\|\bDelta_k\|_2^2])^{1/2} + (\bfE[\|\bV_k\|_2^2])^{1/2}\\
		&\le \gamma (\bfE[\|\bDelta_k\|_2^2])^{1/2} + 1.82(h^3 M^2 p)^{1/2},
\end{align}
where we have used the fact that $h\le 2/M$.
Using this inequality iteratively with $k-1,\ldots,0$ instead of $k$, we get
\begin{align}
(\bfE[\|\bDelta_{k+1} \|_2^2])^{1/2}  
		&\le \gamma^{k+1} (\bfE[\|\bDelta_0\|_2^2])^{1/2} + 1.82 (h^3 M^2 p)^{1/2}
		\sum_{j=0}^{k} \gamma^j\\
		&\le \gamma^{k+1} (\bfE[\|\bDelta_0\|_2^2])^{1/2} + 1.82 (h^3 M^2 p)^{1/2}
		(1-\gamma)^{-1}.\label{F}
\end{align}
Since $\bDelta_{k+1}  = \bL_{(k+1)h} - \bvartheta^{(k+1,h)}$ and $\bL_{(k+1)h}\sim \pi$, 
we readily get the inequality $W_2(\nu_{k+1},\pi)\le \big(\bfE[\|\bDelta_{k+1} \|^2_2]\big)^{1/2}$. 
In addition, one can choose $\bL_0$ so that $W_2(\nu_0,\pi)=\big(\bfE[\|\bDelta_{0} \|^2_2]\big)^{1/2}$.
Using these relations and substituting $\gamma$ by its expression in \eqref{F},
we get the two claims of the theorem.

\subsection{Proof of \Cref{thTwo}}
Using \eqref{C}, \eqref{D} and \Cref{lemA}, we get (for every $t>0$)
\begin{align}
\bfE[\|\bDelta_{k+1} \|_2^2]
		&= \bfE[\|\bDelta_k  -h \bU_k + \bV_k\|_2^2] + \bfE[\|\sigma h \bzeta^{(k)}\|_2^2]\\
		&\le (1+t)\bfE[\|\bDelta_k  -h \bU_k\|_2^2] + (1+t^{-1})\bfE[\|\bV_k\|_2^2]+ \sigma^2 h^2 p \\
		&\le (1+t)\gamma^2\bfE[\|\bDelta_k\|_2^2] + (1+t^{-1})\bfE[\|\bV_k\|_2^2]+ \sigma^2 h^2 p  .
\end{align}
Since $h\le 2/M$, \Cref{lemC} implies that 
\begin{align}
\bfE[\|\bDelta_{k+1} \|_2^2]
		&\le (1+t)\gamma^2\bfE[\|\bDelta_k\|_2^2] +   (1+t^{-1})(1.82)^2h^3M^2p + \sigma^2 h^2 p
\end{align}
for every $t>0$. Let us choose $t = (\frac{1+\gamma}{2\gamma})^2 - 1$ so that 
$(1+t)\gamma^2 = (\frac{1+\gamma}{2})^2$. By recursion, this leads to
\begin{align}
W_2^2(\nu_{k+1},\pi)
		&\le \Big(\frac{1+\gamma}{2}\Big)^{2(k+1)} W_2^2(\nu_0,\pi) + 
		\Big(\frac{2}{1-\gamma}\Big)\Big\{\sigma^2 h^2 p + (1+t^{-1})(1.82)^2h^3M^2p\Big\}.
\end{align}
In the case $h\le 2/(m+M)$,  $\gamma = 1-mh$ and  we get $\frac{1+\gamma}{2} = 1-\frac12 mh$. Furthermore,
\begin{align}
(1+t^{-1})h^3M^2p & = \frac{(1+\gamma)^2h^3M^2p}{(1-\gamma)(1+3\gamma)}\le 
\frac{h^2M^2p}{m}.
\end{align}  
This readily yields
\begin{align}
W_2(\nu_{k+1},\pi)
		&\le \Big(1-\frac{mh}{2}\Big)^{k+1} W_2(\nu_0,\pi) + 
		\Big(\frac{2hp}{m}\Big)^{1/2}\Big\{\sigma^2  + \frac{3.3M^2}{m}\Big\}^{1/2}.
\end{align}
Similarly, in the case $h\ge 2/(m+M)$,  $\gamma = Mh-1$ and we get 
$\frac{1+\gamma}{2} = \frac12 Mh$. Furthermore,
\begin{align}
(1+t^{-1})h^3M^2p & = \frac{(1+\gamma)^2h^3M^2p}{(1-\gamma)(1+3\gamma)}\le 
\frac{h^3M^2p}{2-Mh} \le \frac{2h^2Mp}{2-Mh}.
\end{align}  
This implies the inequality
\begin{align}
W_2(\nu_{k+1},\pi)
		&\le \Big(\frac{Mh}{2}\Big)^{k+1} W_2(\nu_0,\pi) + 
		\Big(\frac{2h^2p}{2-Mh}\Big)^{1/2}\Big\{\sigma^2  + \frac{6.6M}{2-Mh}\Big\}^{1/2},
\end{align}
which completes the proof. 

\subsection{Proofs of lemmas}
\label{ssecLem}

\begin{proof}[Proof of \Cref{lemA}]
Since $f$ is $m$-strongly convex, it satisfies the inequality
\begin{equation}
\bDelta^\top \big(\nabla f(\bvartheta+\bDelta) - \nabla f(\bvartheta)\big)
\ge \frac{mM}{m+M}\|\bDelta\|_2^2 + \frac1{m+M} \|\nabla f(\bvartheta+\bDelta) - \nabla f(\bvartheta)\|_2^2,
\end{equation}
for all $\bDelta,\bvartheta\in\RR^p$. Therefore, simple algebra yields
\begin{align}
\|\bDelta_k  -h \bU_k\|_2^2 
		& = \|\bDelta_k\|_2^2 - 2h \bDelta_k^\top \bU_k+ h^2 \|\bU_k\|_2^2 \\
		& = \|\bDelta_k\|_2^2 - 2h \bDelta_k^\top \big(\nabla f(\bvartheta^{(k,h)}+\bDelta_k)
				-\nabla f(\bvartheta^{(k,h)})\big) + h^2 \|\bU_k\|_2^2 \\
		& \le \|\bDelta_k\|_2^2 - \frac{2h mM}{m+M} \|\bDelta_k\|_2^2 - 
				\frac{2h}{m+M} \|\bU_k\|_2^2 + h^2 \|\bU_k\|_2^2 \\			
		& = \Big(1 - \frac{2h mM}{m+M}\Big) \|\bDelta_k\|_2^2 + 
				h\Big(h - \frac{2}{m+M}\Big)\|\bU_k\|_2^2.\label{G}
\end{align}
Note that, thanks to the strong convexity of $f$, the inequality 
$\|\bU_k\|_2 = \|\nabla f(\bvartheta^{(k,h)}+\bDelta_k) - \nabla f(\bvartheta^{(k,h)}) \|_2\ge 
m\|\bDelta_k\|_2$ is true. If $h\le \nicefrac2{(m+M)}$, this inequality can be
combined with \eqref{G} to obtain
\begin{align}
\|\bDelta_k  -h \bU_k\|_2^2 
		& \le  (1 - h m)^2 \|\bDelta_k\|_2^2
		.\label{H}
\end{align}
Similarly, when $h\ge \nicefrac2{(m+M)}$,  we can use the Lipschitz property of $\nabla f$
to infer that $\|\bU_k\|_2 \le M\|\bDelta_k\|_2$. Combining with \eqref{G}, this
yields
\begin{align}
\|\bDelta_k  -h \bU_k\|_2^2 
		& \le  (h M-1)^2 \|\bDelta_k\|_2^2,\qquad \text{if}\qquad h\ge \nicefrac2{(m+M)}.\label{I}
\end{align}
Thus, we have checked that \eqref{E} is true for every $h\in(0,\nicefrac2M)$.
\end{proof}

\begin{proof}[Proof of \Cref{lemB}]
To simplify notations, we prove the lemma for $p=1$. The function 
$x\mapsto f'(x)$ being Lipschitz continuous is almost surely differentiable. 
Furthermore, it is clear that $|f''(\bx)|\le M$  for every $\bx$ for which 
this second derivative exists. The result of \citep[Theorem 7.20]{rudin87} 
implies that 
\begin{equation}
f'(x)-f'(0) = \int_0^x f''(y)\,dy.
\end{equation}  
Therefore, using $f'(x)\,\pi(x) = -\pi'(x)$, we get
\begin{align}
\int_{\RR} f'(x)^2\,\pi(x)\,dx 
		& = f'(0)\int_{\RR} f'(x)\,\pi(x)\,dx + 
				\int_{\RR}\Big(\int_0^x f''(y)\,dy\Big) f'(x)\,\pi(x)\,dx \\
		& = -f'(0)\int_{\RR} \pi'(x)\,dx - 
				\int_{\RR}\Big(\int_0^x f''(y)\,dy\Big) \pi'(x)\,dx \\
		& = -\int_{0}^\infty\int_0^x f''(y)\,\pi'(x)\,dy\,dx
				+\int_{-\infty}^0\int_x^0 f''(y)\,\pi'(x)\,dy\,dx.
\end{align}
In view of Fubini's theorem, we arrive at
\begin{align}
\int_{\RR} f'(x)^2\,\pi(x)\,dx & =
		\int_{0}^\infty f''(y)\,\pi(y)\,dy
				+\int_{-\infty}^0 f''(y)\,\pi(y)\,dy\le M.
\end{align}
This completes the proof.
\end{proof}

\begin{proof}[Proof of \Cref{lemC}] 
Since the process $\bL$ is stationary, $V(a)$ has the same distribution as $V(0)$. For this
reason, it suffices to prove the claim of the lemma for $a=0$ only.
Using the Lipschitz continuity of $f$, we get
\begin{align}
\bfE[\|\bV(0)\|^2_2]
		& = \bfE\Big[\Big\|\int_{0}^{h}\big(\nabla f(\bL_t) - \nabla f(\bL_{0})\big)\,dt\Big\|^2_2\Big]\\
		& \le  h\int_{0}^{h}\bfE\big[\big\|\nabla f(\bL_t) - \nabla f(\bL_{0})\big\|^2_2\big]\,dt\\
		& \le  hM^2\int_{0}^{h}\bfE\big[\big\|\bL_t - \bL_{0}\big\|^2_2\big]\,dt.
\end{align}		
Combining this inequality with the stationarity of $\bL_t$, we arrive at
\begin{align}		
\Big(\bfE[\|\bV(0)\|^2_2]\Big)^{1/2}
		& \le \bigg(hM^2\int_{0}^{h}\bfE\big[\big\|-\int_{0}^t \nabla f(\bL_s)\,ds + 
				\sqrt{2}\,\bW_{t}\big\|^2_2\big]\,dt\bigg)^{1/2}\\
		& \le \bigg(hM^2\int_{0}^h\bfE\big[\big\|\int_{0}^t \nabla f(\bL_s)\,ds\big\|^2_2\big]\,dt\bigg)^{1/2}
				+ \bigg(2hpM^2\int_0^h t\,dt\bigg)^{1/2}\\
		& \le \bigg(hM^2\bfE\big[\big\|\nabla f(\bL_0)\big\|^2_2\big]\int_0^h t^2\,dt\bigg)^{1/2}
				+ \bigg(2hpM^2\int_0^h t\,dt\bigg)^{1/2}\\				
    & = \bigg(\frac13 h^4M^2 \bfE\big[\big\|\nabla f(\bL_0)\big\|^2_2\big]\bigg)^{1/2}
				+ \big(h^3 M^2p\big)^{1/2}. 
\end{align}
To complete the proof, it suffices to apply \Cref{lemB}.
\end{proof}

{\renewcommand{\addtocontents}[2]{}

\acks{
The work of the author was partially supported by the grant 
Investissements d'Avenir (ANR-11-IDEX-0003/Labex Ecodec/ANR-11-LABX-0047). The author would 
like to thank Nicolas Brosse, who suggested an improvement in \Cref{thTwo}. 
}%

{\renewcommand{\addtocontents}[2]{}
\bibliography{Literature}}

\end{document}